\documentclass[10pt]{amsart}

\usepackage{latexsym, amsmath, amsthm, amsopn, amsfonts,graphicx, amssymb,epsfig,graphics,hyperref}
\usepackage{tikz}

\newtheorem{theo}{Theorem}[section]

\newtheorem{lemma}[theo]{Lemma}

\newtheorem{rem}[theo]{Remark}

\newcommand{\R}{\mathbb{R}}
\newtheorem{thm}{Theorem}[section]

\begin{document}
\title[{Quantitative isoperimetric inequality with barycentric distance}]{Optimal sets for the quantitative isoperimetric inequality in the plane with the barycentric distance}
               
\author{Gisella Croce}
\author{Antoine Henrot} 

\address{{G. Croce:} 
SAMM, UR 4543 Universit\'e Paris 1 Panth\'eon-Sorbonne, FR 2036 CNRS, F-75013 Paris, France }
\email{gisella.croce@univ-paris1.fr}

\address{{A. Henrot:} 
Universit\'e de Lorraine, CNRS, IECL, F-54000 Nancy,  France }
\email{antoine.henrot@univ-lorraine.fr}


\date{\today}


\begin{abstract}
In a recent paper, C. Gambicchia and A. Pratelli proved a quantitative isoperimetric inequality involving the 
isoperimetric deficit $\delta(K)$ and the barycentric distance $\lambda_0(K)$
for sets $K\subset \R^N$ with given diameter $D$ and measure. In this work we are interested in the optimal sets for this inequality in the plane,
i.e. sets that minimize the ratio $\delta(K)/\lambda_0(K)^2$. We prove exi\-sten\-ce of optimal sets (at least when $D$ is large enough),
 regularity and express the optimality conditions. Moreover, we prove that the optimal sets have exactly two connected
 components and their boundary does not contain any arc of circle.

\end{abstract}
\maketitle

\section{Introduction}

Quantitative isoperimetric inequalities have received much attention in the last thirty years literature.
The question is to quantify the distance that a subset of $\R^n$, $\Omega$, has from an $n$-dimensional ball $B$ of the same measure in terms of the isoperimetric deficit 
\begin{equation}\label{def-delta}
\delta(\Omega)=\frac{P(\Omega)-P(B)}{P(B)}\,
\end{equation}
through  inequalities of the type 
\begin{equation}\label{distK}
[dist(\Omega,B)]^k\le C\; \delta(\Omega)\,.
\end{equation}
 In terms of distances, basically, the Hausdorff distance and the measure of the symmetric difference between 
 $\Omega$ and the ball of same measure have been proposed in the litterature. Notice that in both case, on can consider either 
 the ball centered at the barycentre $x^G$ of $\Omega$, or the ball minimizing the distance of $\Omega$ from all the balls of $\R^n$ with fixed volume.

More in details, the Hausdorff distance has been considered in  \cite{Fu89Transactions} and later in \cite{FGP}. However, it requires some geometrical structure on the sets (otherwise one can consider  a "big" ball and a very tiny ball: moving away the two balls from each other, the isoperimetric deficit  
can be made arbitrarily small, and at the same time  the Hausdorff asymmetry  can be made arbitrarily large). 
To study more general families of sets,  one can use $L^1$ distances between sets.
L. E. Fraenkel proposed the now so called Fraenkel asymmetry 
\begin{equation}\label{def-lambda}
\lambda(\Omega)=\inf_{y\in \R^n}\frac{|\Omega \Delta B_y|}{|\Omega|}\,
\end{equation}
where $B_y$ is the ball of center $y$,  having the same Lebesgue measure  of  
$\Omega$ and $\Delta$ is the sym\-me\-tric difference of sets.
Many mathematicians studied quantitative isoperimetric inequalities with the Fraenkel asymmetry, establishing sharp inequalities (see for example \cite{HHW},  \cite{H}, \cite{Ca}, \cite{AFN}, \cite{FiMP}, \cite{FuMP}, \cite{CiLe}, \cite{Fuscopreprint}, {\cite{DL}})
 and existence of an optimal set for the related shape optimization problem (see \cite{CiLeexistence} and \cite{BCH_COCV}).

Fuglede  proposed in \cite{Fu93Geometriae} 
the barycentric asymmetry:
$$
\lambda_0(\Omega)=\frac{|\Omega \Delta B_{x^G}|}{|\Omega|}\,,
$$
where $B_{x^G}$ is the ball centered at the barycentre ${x^G}$ of $\Omega$ and such that $|\Omega|=|B_{x^G}|$.  
Notice that $\lambda_0(\Omega)$ is obviously
easier to compute than $\lambda(\Omega)$, since it does not require a minimization problem to solve.
Fuglede proved that there exists a positive constant (depending only on the dimension $n$) such that
\begin{equation}\label{Fuglede_convex}
\delta(\Omega)\geq C(n)\,\lambda_0^2(\Omega),\quad\text{for any convex subsets $\Omega$ of $\R^n$}.
\end{equation}
This inequality holds more in general for compact connected  sets of the plane, as  proved
in \cite{BCH_Annali}. We notice that  one needs to add a constraint, such as the connectedness of the sets,   to have a well-posed problem. 
Otherwise, the union of a unit ball and a second ball with radius $r<<1$ 
at distance $d>>1$ has a very small isoperimetric deficit while the barycentric asymmetry can be kept equal to 2.
A different, very natural constraint to add, in order to be able to establish a meaningful inequality, is the sets to be  contained in a fixed  ball.
In \cite{GP24},
the authors proved the following 
\begin{thm}
For every
$N\geq 2$  and every $D>0$
there exists a constant $C(N,D)$ 
such that, 
for any set 
$E\subset \R^N$
with
diameter less than $D |E|^{1/N}$, the inequality
$\lambda^2_0(E)\leq  C(N, D)\delta(E)$ 
holds true.
\end{thm}
Note that instead of bounding the product $D |E|^{1/N}$,we can choose to fix the area and bound the diameter separately.
This result gave us the idea of investigating the corresponding minimization problem
$$
\min\{J(K), K \mbox{ compact},  K\subset\mathbb{R}^2, diam(K)\leq D, |K|=\pi\}
$$
where
\begin{equation}\label{defJ}
J(K)=\frac{\delta(K)}{\lambda^2_0(K)}
\end{equation}
studying the existence and  properties of the minimizer. We gather here 
our main results:
\begin{thm}
\begin{enumerate}
\item
There exists a number $D^*$, with $D^* \leq 10$ such that, for any $D\geq D^*$, the  above minimization problem 
has a solution.
\item
The boundary of the minimizer  is globally $C^{1,1}$;
the parts of the boundary that do not touch the barycentric ball are analytic.
\item
There exists a number $D^{**}$ such that for $D\geq D^{**}$, the minimizer has exactly two connected components,  touching the boundary of the ball of diameter $D$
in opposite points.
\item
The boundary of the minimizer does not contain any arc of circle. 
\end{enumerate}
\end{thm}
This last property  is quite surprising, if one thinks that for example the minimizer of 
$
 \frac{\delta(\Omega)}{\lambda^2(\Omega)}
$
is conjectured to be a "mask",  whose boundary is composed by arcs of circle, with three different curvatures (see \cite{CiLeexistence, BCH_COCV}). In the case of convex sets, the minimizer is an explicitly described stadium, as proved in \cite{AFN}.
We will also derive a formula of the curvature of the optimal set (see Theorem \ref{ocPrinceton} for all the details) .

The paper is organized as follows.
After some preliminary results,  we  prove in section 3 the existence of a minimizing set of $J(K)$
for all compact $K\subset \R^2$ contained in a sufficiently large ball.
 In section \ref{section4} we will study the regularity of the minimizer and write an optimality condition involving its curvature using the shape derivative. 
Section 5 is devoted to the number of connected components of the optimal set and the proof  that its boundary  does not contain any arc of circle.
\section{Preliminaries}\label{section_preliminaries}
We recall in this section some definitions and results that will be useful in the whole paper.

The barycenter of a set $\Omega$ is  defined as
$\displaystyle
x^G=\frac{1}{|\Omega|}\int_{\Omega} x\,dx\,.
$
For a set $\Omega\subset \R^n$,  $|\Omega|$ denotes the 
 $n$-dimensional Lebesgue measure. 

As proved in \cite{CiLeexistence} and  \cite{BCH_COCV} for a sequence $K_\varepsilon$ converging to a ball, one has
\begin{theo}\label{thmBCH}
Let $\{K_{\varepsilon} \}_{\varepsilon>0}$ be a sequence of compact planar sets converging
to a ball $B$ in the sense that  
$|B\Delta K_{\varepsilon} |\to 0$ as $\varepsilon\to 0$.
Then
$$
\inf
\left\{
\liminf_{\varepsilon\to 0}\frac{\delta(K_{\varepsilon})}{\lambda^2(K_{\varepsilon})}
\right\} 
= 
\frac{\pi}{8(4-\pi)}\simeq 0.457474
$$
where $\lambda$ is the Fraenkel asymmetry, defined in (\ref{def-lambda}).
\end{theo}
The following  lemma, proved in \cite{BCH_Annali},  gives the $L^1$ distance of two balls in terms of the distance of their centers:
 \begin{lemma}\label{lemmapalle}
Let $B_{(a,0)}$ denote the ball of area $\pi$ centered at $(a,0)$, $a>0$.  We have 
$$
|B_{(0,0)}\Delta B_{(0,a)}|=d_{L^1}(B_{(0,0)},B_{(0,a)})=4\arcsin\left(\frac{a}{2}\right)+2a\sqrt{1-\frac{a^2}{4}}=4a + o(a).
$$
\end{lemma}
In the introduction we have already explained why the assumption of the sets being in ball of diameter $D$ makes the minimisation of $J$ meaningful. We make the example precise and use the conclusion of the following
\begin{rem}\label{limite_Fuglede}
Let 
$\Omega_n$ be the union of 
the disk centered in $(2,0)$, of radius $R_n=1-\frac{1}{n}$, and the disk
centered in
$\left(-\frac{2(n-1)^2}{2n-1},0\right)$,
 of radius
$r_n=\sqrt{\frac{2n-1}{n^2}}$\,.
It is easy to check that $|\Omega_n|=\pi$, the barycentre of $\Omega_n$ is the origin, $\delta(\Omega_n)=R_n+r_n-1\to 0$ as $n\to \infty$
and $\lambda_0(\Omega_n)=2$. Thus $\lim_{n\to\infty} J(\Omega_n)=0$.
In other words,
$$\displaystyle \lim_{D\to +\infty} \inf\{J(K), diam(K)\leq D\} =0.
$$
\end{rem}
We will use 
that, in dimension two, for a connected set, the 
perimeter of the convex hull is smaller than the
perimeter of the set itself.  This is proved for example in \cite{Ferriero-Fusco}.

%
%
%

We are now going to prove  a Lemma that can be found in \cite{GP24}.
We state it in a way that will be convenient for us.
\begin{lemma}\label{lemgrav}
Let $E,F$ be two compact sets with the same measure. Let $G_E,G_F$ be their barycenter, then
$$|G_E-G_F|\leq \frac{D'}{2|E|} |E \Delta F|$$
where $D'$ is the diameter of $E \cup F$.
\end{lemma}
\begin{proof}
The proof  is immediate, observing that
\begin{eqnarray*}
|G_E-G_F| =\frac{1}{|E|}\left|\int_E x dx - \int_F x dx\right| = \frac{1}{|E|}\left|\int_{E\setminus F} x dx - \int_{F\setminus E} x dx\right| \\
\leq \frac{D'}{|E|} |E\setminus F| = \frac{D'}{2|E|} |E \Delta F|.
\end{eqnarray*}
\end{proof}
To study the regularity of the minimizer we will use the notion of
{\it strong $\Lambda$-minimizer of the perimeter} (that is similar to the classical notion of quasi-minimizer of the perimeter).
We say that a set $E$ is a strong $\Lambda$-minimizer of the perimeter if there exists $R>0$ such that for any $x\in \mathbb{R}^2$,
for any $r, 0<r \leq R$ and for any set $F$ such that $E\Delta F \subset\subset B(x,r)$ (i.e. $F$ is a compact variation of $E$ inside the
ball centered at $x$ of radius $r$), we have
\begin{equation}
P(E\cap B(x,r)) \leq P(F\cap B(x,r)) + \Lambda |E\Delta F|.
\end{equation}
It is known that a strong $\Lambda$-minimizer of the perimeter is $C^{1,1}$ in dimension two, see e.g.  F. Maggi's book \cite{Maggi} and also \cite{CiLeexistence}.
We  will use Theorem 3.6 of \cite{stred-zi} and  classical results for the
isoperimetric problem
(see for example \cite{tamanini}):
\begin{thm}\label{thm_stred-zi}
Assume that $\Omega \subset \R^2$ is bounded, convex and has a $C^2$ boundary. 
If $F$ is a solution
of
$
\inf\{P(E,\R^2), E\subset \overline{\Omega}, |E|=v < |\Omega|\}
$,
then $F$ is $C^{1,1}$ in some neighborhood of 
$\partial \Omega$.
Moreover $\partial F\cap \Omega$
is real analytic.
\end{thm}

\section{Existence}
In this section, we prove existence of a minimizer for the functional $J(K)=\delta(K)/\lambda_0^2(K)$ in the class of compact sets
of area $\pi$ and diameter less than $D$, at least when $D$ is large enough.  More precisely:
\begin{theo}\label{theoexistence}
There exists a number $D^*$, with $D^* \leq 10$ such that, for any $D\geq D^*$, the following minimization problem has a solution
\begin{equation}
\min\{\delta(K)/\lambda_0^2(K), K \mbox{ compact},  K\subset\mathbb{R}^2, diam(K)\leq D, |K|=\pi\}.
\end{equation}
\end{theo}
The main difficulty is to exclude a minimizing sequence $K_\varepsilon$ converging to the disk. For that purpose, we will give an estimate (from below) of $\delta(K_\varepsilon)/\lambda_0^2(K_\varepsilon)$: we will prove that
 $$
 \liminf \frac{\delta(K_\varepsilon)}{\lambda_0^2(K_\varepsilon)} \simeq 0.0885
 $$ 
 (independently of the diameter $D$). Moreover we will explicitly describe a set $K$
such that $J(K)$ is lower than 0.0885. $K$ is the union of two disks, the biggest disk being tangential to the barycentric disk. At this step we will need a sufficiently
large diameter $D$ to get a value of $J$ lower than the bound $0.0885$ found before.

\begin{proof}[Proof of Theorem \ref{theoexistence}]
In the sequel, we assume that the diameter $D$ satisfies $D\geq 5 > 2+2\sqrt{2}$. Only in the last step of the proof we will need $D\geq 10$.\\

Let $K_\varepsilon$ be a minimizing sequence, that is,  $\displaystyle \frac{\delta(K_\varepsilon)}{\lambda_0^2(K_\varepsilon)} \to \inf_E J(E) $.
Without loss of generality, we can assume that all the sets $K_\varepsilon$ have area $\pi$.
Since $D> 2+2\sqrt{2}$, an admissible set $K$ is composed of two balls of same radius $1/\sqrt{2}$ with their centers at a distance $2+\sqrt{2}$.
In that case the barycentric ball is tangent to the two balls and $\lambda_0(K)=2$. Since $\delta(K)=\sqrt{2}-1$ we see that $J(K)=(\sqrt{2}-1)/4$ and then we can assume that 
\begin{equation}\label{firstestim}
J(K_\varepsilon) \leq \frac{\sqrt{2}-1}{4} .
\end{equation}
Since $\lambda_0(E)\leq 2$ for any set $E$ and $\delta(K_\varepsilon)=\frac{P(K_\varepsilon)}{2\pi} -1$,  we get from (\ref{firstestim})
\begin{equation}\label{perimeterbounded}
P(K_\varepsilon)\leq 2\pi \sqrt{2} \leq 9.
\end{equation} 
We will consider the case where the minimizing sequence converges to a compact set $K$ distinct of a ball at the end of the proof
(and this limit set will provide the existence of a minimizer). Therefore, let us concentrate on the case where $K_\varepsilon$ converges
to a ball.

\begin{enumerate}
\item[Step 1]
We can assume that
$\delta(K_\varepsilon)\to 0$ and $\lambda(K_\varepsilon)=2\varepsilon\to 0$.
By Theorem \ref{thmBCH} one has
$$
\delta(K_\varepsilon)\geq 0.4574 \cdot 4 \varepsilon^2\,.
$$
We are now going to prove that
\begin{equation}\label{difflambda}
\lambda_0(K_\varepsilon)-\lambda(K_\varepsilon)= |\lambda_0(K_\varepsilon)-\lambda(K_\varepsilon)|\leq \frac{4A}{\pi}\varepsilon
\end{equation}
for some explicit constant $A>0$.
This will imply that
\begin{equation}
\frac{\delta(K_\varepsilon)}{\lambda_0^2(K_\varepsilon)}
\geq
\frac{\delta(K_\varepsilon)}{\left(\lambda(K_\varepsilon)+\frac{4A}{\pi}\varepsilon\right)^2}
\geq \frac{1.8296}{\left(2+\frac{4A}{\pi}\right)^2}\,.
\end{equation}
This gives an useful estimate in the case of a minimizing sequence $K_\varepsilon$ converging to a ball that we will use later.

To prove (\ref{difflambda})  it is sufficient to find  a positive constant $A$ such that 
\begin{equation}
\label{distance-centres-a-prouver}
|G_\varepsilon - F_\varepsilon| \leq A \varepsilon
\end{equation}
where $G_\varepsilon$ is the barycentre of $K_\varepsilon$ and $F_\varepsilon$ is the centre of an optimal (Fraenkel) ball for $\lambda(K_\varepsilon)$.
Indeed, by the triangle inequality, 
$$
\pi(\lambda_0(K_\varepsilon)-\lambda(K_{\varepsilon}))=d_{L^1}(K_\varepsilon,B_{G_\varepsilon})- d_{L^1}(K_\varepsilon, B_{F_\varepsilon})\leq d_{L^1}(B_{G_\varepsilon}, B_{F_\varepsilon})\,,
$$
where $B_{G_\varepsilon}$ is the barycentric ball and $B_{F_\varepsilon}$ is an optimal ball for the Fraenkel asymmetry, both for $K_\varepsilon$.
This inequality together with  (\ref{distance-centres-a-prouver}) and Lemma \ref{lemmapalle} imply (\ref{difflambda}).

Using Lemma \ref{lemgrav} with $E=K_\varepsilon$ and $F=B_{F_\varepsilon}$ is its Fraenkel ball, we obtain \eqref{distance-centres-a-prouver} with $A=D$
(note that $\lambda(K_\varepsilon)=2\varepsilon=|K_\varepsilon \Delta B_{F_\varepsilon}|/\pi$).
This gives the first estimate
\begin{equation}\label{roughestimate}
\frac{\delta(K_\varepsilon)}{\lambda_0^2(K_\varepsilon)}
\geq
 \frac{1.8296}{\left(2+\frac{4D}{\pi}\right)^2}\,.
\end{equation}
In the next steps we are going to estimate $D$.

To sum up what we obtain at the previous step, the sequence $K_\varepsilon$ (converging to a ball) satisfies
$$2\varepsilon \leq \lambda_0(K_\varepsilon) \leq \left(2+\frac{4D}{\pi}\right) \varepsilon$$
by (\ref{difflambda}), (\ref{distance-centres-a-prouver}) with $A=D$ and $\lambda(K_{\varepsilon})=2\varepsilon$.

This and estimate \eqref{firstestim} give $\delta(K_\varepsilon) \leq \frac{\sqrt{2} -1}{4} \left(2+\frac{4D}{\pi}\right)^2 \varepsilon^2.$
Since $\lambda(K_{\varepsilon})=2\varepsilon$,
by Theorem \ref{thmBCH} we get $\delta(K_{\varepsilon})\geq 0.4574\cdot(2\varepsilon)^2$ and therefore
\begin{equation}\label{estidelta}
1.8296 \varepsilon^2 \leq \delta(K_\varepsilon) \leq \frac{\sqrt{2} -1}{4} \left(2+\frac{4D}{\pi}\right)^2 \varepsilon^2.
\end{equation}
In particular, we see that $\delta(K_\varepsilon)=O(\varepsilon^2)$.
Without loss of generality, we will assume now that the center of the Fraenkel ball $B_{F_\varepsilon}$ remains at the origin
(up to some translation and rotation). Then, we will denote by $B_F$ the Fraenkel ball that remains fixed.

\item[Step 2]
We will study the connected components of $K_\varepsilon$.
Let us consider a connected component of $K_\varepsilon$, say $K_\varepsilon^1$ of measure $\pi m_\varepsilon$ and let us denote by $K^0_\varepsilon$
its complement in $K_\varepsilon$ of measure $\pi (1-m_\varepsilon)$.
By the classical isoperimetric inequality applied to these two sets, we have
$$P^2(K^1_\varepsilon) \geq 4\pi^2 m_\varepsilon \quad\mbox{ and }  \quad P^2(K^0_\varepsilon) \geq 4\pi^2 (1-m_\varepsilon).$$
Therefore
$$\delta(K_\varepsilon)=\frac{P(K^0_\varepsilon)+P(K^1_\varepsilon)-2\pi}{2\pi} \geq \sqrt{1-m_\varepsilon} +  \sqrt{m_\varepsilon} -1.$$
From the last inequality in 
\eqref{estidelta} we infer that $m_\varepsilon$ or $1-m_\varepsilon$ must be of order $\varepsilon^4$:
\begin{equation}\label{ordere4}
m_\varepsilon=O(\varepsilon^4) \quad\mbox{ or }  \quad 1-m_\varepsilon=O(\varepsilon^4) .
\end{equation}
There is exactly one "big" connected component that meets the Fraenkel 
ball $B_F$ and its measure is, according to
\eqref{ordere4}, $\pi-O(\varepsilon^4)$. Let us denote by $K^0$ this connected component (for simplicity, we remove the index $\varepsilon$
in what follows) and we are going to prove that the
diameter of $K^0 \cup B_F$ is of order $2+O(\varepsilon^{2/3})$. 

\item[Step 3]
We will give an estimate of the diameter of $K^0$.
We recall that $|B_F\setminus K_\varepsilon|=\pi \varepsilon$, as $\lambda(K_{\varepsilon})=2\varepsilon$. Decomposing $K_\varepsilon=K^0 \cup K^1$ with $|K^1|=O(\varepsilon^4)$,
we see that 
$$B_F\setminus K^0  = (B_F\setminus K_\varepsilon) \cup (B_F\cap K^1)/,.$$
Therefore
\begin{equation}
\pi\varepsilon \leq |B_F \setminus K^0| \leq \pi\varepsilon + O(\varepsilon^4) < 4\varepsilon .
\end{equation}
This estimate implies that in each subdomain of $B_F$ with measure greater than $4\varepsilon$ we can find at least one point of $K^0$.
Let us consider for example a domain limited by the boundary of the ball $B_F$ and a chord whose aperture is $2\alpha_\varepsilon$.
We call it a {\it $\varepsilon$-cap} in the sequel.
The area of such a domain is $A_\varepsilon=\alpha_\varepsilon - \sin \alpha_\varepsilon \cos \alpha_\varepsilon$ and for $\varepsilon$ small
we have $A_\varepsilon = 2\alpha_\varepsilon^3/3 + o(\varepsilon^3)$. Thus, by choosing $\alpha_\varepsilon =(6\varepsilon)^{1/3}$ we know that
each $\varepsilon$-cap  contains a point of $K^0$. This implies that the convex hull of $K^0$, containing a point in each possible $\varepsilon$-cap, 
certainly contains the disk tangent to the chords
joining the extremities of these domains (see Figure \ref{figure1}),
namely the disk $D_\varepsilon$ centered at the origin of radius $R_\varepsilon:=\cos(2\alpha_\varepsilon)$ that satisfies, by Taylor expansion
\begin{equation}\label{estimateReps}
R_\varepsilon = 1-2\alpha_\varepsilon^2 +o(\varepsilon^2)=1-2(6\varepsilon)^{2/3} +o(\varepsilon^{2/3}).
\end{equation}
\begin{center}
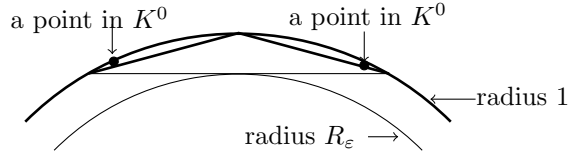
\begin{figure}[h]
\begin{tikzpicture}
    \draw [thick, line width=1pt] (6,3) arc [radius=4cm, start angle=45, end angle= 135];
    \draw [thick, line width=1pt] (5.18,3.64) -- (3.18,4.18)--(1.18,3.64);
    \draw (5.18,3.64) -- (1.18,3.64);
     \draw (5.62,2.62) arc [radius=3.46cm, start angle=45, end angle= 135];
     \draw [<-] (5.73,3.3) -- (6.33,3.3);
     \node at (6.96,3.34) {radius $1$};
      \draw [->] (4.9,2.8) -- (5.3,2.8);
     \node at (4,2.8) {radius $R_\varepsilon$};
     \draw [fill] (4.85,3.75) circle [radius=0.06];
      \node at (4.9,4.4) {a point in $K^0$};
     \draw [->] (4.85,4.2)--(4.85,3.8);
     \draw [fill] (1.52,3.8) circle [radius=0.06];
     \node at (1.2,4.35) {a point in $K^0$};
     \draw [->] (1.5,4.3)--(1.5,3.9);
\end{tikzpicture}
\caption{The convex hull of $K^0$ contains the disk of radius $R_\varepsilon=1-2(6\varepsilon)^{2/3}$}\label{figure1}
\end{figure}
\end{center}
Now, let $M$ be a point of $K^0$ of maximal distance from the origin. By definition the convex hull $conv(K^0)$ of $K^0$ contains the convex hull of
$M  \cup D_\varepsilon$, say $conv(M,D_\varepsilon)$.  Using the relation between the perimeters recalled above
$$
P(K^0) \geq P(conv(K^0)) \geq P(conv(M,D_\varepsilon)).
$$
This last perimeter can be computed explicitly: indeed the perimeter of the convex hull of a disk of radius $R_1$ and a point at a distance $R_2>R_1$
from the center of the disk equals the following functions of $R_1, R_2$:
$$
p(R_1,R_2)= R_1\left(2\pi -2 \arccos\left(\frac{R_1}{R_2}\right)\right) +2\sqrt{R_2^2-R_1^2}.
$$
In our case,  by (\ref{estimateReps}), $R_1=1-2(6\varepsilon)^{2/3}+o(\varepsilon^{2/3})$, and we can write, by (\ref{estidelta})
\begin{equation}\label{perimetro}
\frac{\sqrt{2} -1}{4} \left(2+\frac{4D}{\pi}\right)^2 \varepsilon^2 \geq \delta(K_\varepsilon) \geq 
\frac{p(R_1,R_2)}{2\pi} -1\,.
\end{equation}
A Taylor expansion of the function
$(u,v)\mapsto p(1-u,1+v)$ shows that, at the second order:
$$p(1-u,1+v)=2\pi -2\pi u+\frac{4\sqrt{2}}{3} \left(u+v\right)^{3/2} +o((u+v)^{3/2}).$$
By (\ref{perimetro}), with $u=2(6\varepsilon)^{2/3}$ and $R_2=|OM|$ we have
$$
C\varepsilon^2 \geq 
\frac{2\pi -2\pi 2(6\varepsilon)^{2/3}+\frac{4\sqrt{2}}{3}(2(6\varepsilon)^{2/3}+v)^{3/2}+o((2(6\varepsilon)^{2/3}+v))^{3/2}}{2\pi}-1
$$
and therefore, $v=O(\varepsilon^{4/9})$
and this means, in particular, that the diameter of $K^0$ converges to $2$ when $\varepsilon \to 0$.

\item[Step 4]
We will find a new estimate of the liminf for the sequence $K_\varepsilon$.
We have seen that the diameter of $K^0$ converges to $2$. Moreover $K_\varepsilon=K^0 \cup K^1$ with $K^1$ of measure $O(\varepsilon^4)$.
Let us denote by $G^0$ the barycenter of $K^0$ and by $G^1$ the barycenter of $K^1$.  The barycenter of $K_\varepsilon$ is
$G_\varepsilon=(|K^0| G^0 +|K^1| G^1)/\pi$. Therefore
$$
|G_\varepsilon- O|\leq |G_\varepsilon - G^0| +|G^0 -O|=|G^0-O| +O(\varepsilon^4).
$$

By applying Lemma \ref{lemgrav} with $E=K^0$ and $F$ the ball centered at the origin (that is the center of the ball $B_F$) with the same
measure as $K^0$, i.e. a measure that is $\pi-O(\varepsilon^4)$, we obtain
$$|G^0 -O|\leq \frac{2+O(\varepsilon^{4/9})}{\pi-O(\varepsilon^4)} \frac{|K^0\Delta F|}{2}$$
by step 3. We infer
$$
|G_\varepsilon - O| \leq \left(\frac{1}{\pi} + O(\varepsilon^{4/9}) \right) |K^0 \Delta F|.
$$
Now (we recall that $B_F$ is the Fraenkel ball)
$$ |K^0 \Delta F| \leq |K^0 \Delta K_\varepsilon| + |K_\varepsilon \Delta B_F| +|B_F \Delta F| = 2\pi \varepsilon + O(\varepsilon ^4).$$
Passing to the limit, following the framework of step 1 (in particular the inequalities \eqref{distance-centres-a-prouver} and (\eqref{roughestimate})),
 this provides the following threshold $\tau^*$ as an estimate from below of $\liminf J(K_\varepsilon)$:
\begin{equation}\label{threshold}
\tau^*=\frac{1.8296}{\left(2+\frac{8}{\pi}\right)^2} > 0.0885 .
\end{equation}

\item[Step 5]
We are going to exhibit a  domain $K$ such that $J(K)< \tau^*$, where $\tau^*$ has been defined in (\ref{threshold}): this will show that a minimizing sequence cannot converge
to a ball.
$K$ will be  the union of two balls of 
 radii $R_1 \geq R_2$ such that $R_1^2+R_2^2=1$.
Moreover the two balls will not intersect the barycentric ball: in this way $\lambda_0(K)=2$.

Now,  $\delta(K)=R_1+R_2 -1=R_1+\sqrt{1-R_1^2}-1$ 
and this is a decreasing function for $R_1 \in [1/\sqrt{2},1]$.
Therefore, to have $\delta$ as small as possible for this configuration, the best choice is when $R_1$ is as large as possible, i.e. when 
the barycentric ball is exterior, but tangent to the ball of radius $R_1$ as on Figure \ref{fig2}.
\begin{center}
\begin{figure}[h]
\begin{tikzpicture}
    \draw[dotted,thick] (4,3) circle [radius=4cm];
    \draw[fill=blue,thick, line width=1.5pt] (7,3) circle [radius=1cm] ;
     \draw[fill=blue,thick, line width=1.5pt] (0.5,3) circle [radius=5mm] ;
      \draw (4.8,3) circle [radius=1.2cm] ;
      \node at (4.8,3.2) {barycentric};
      \node at (4.8,2.8) {disk};
      \node at (7.5,5) {ball of diameter $D$};
\end{tikzpicture}
\caption{A competitor obtained with two disks, the barycenter disk being tangent to the big one.}\label{fig2}
\end{figure}
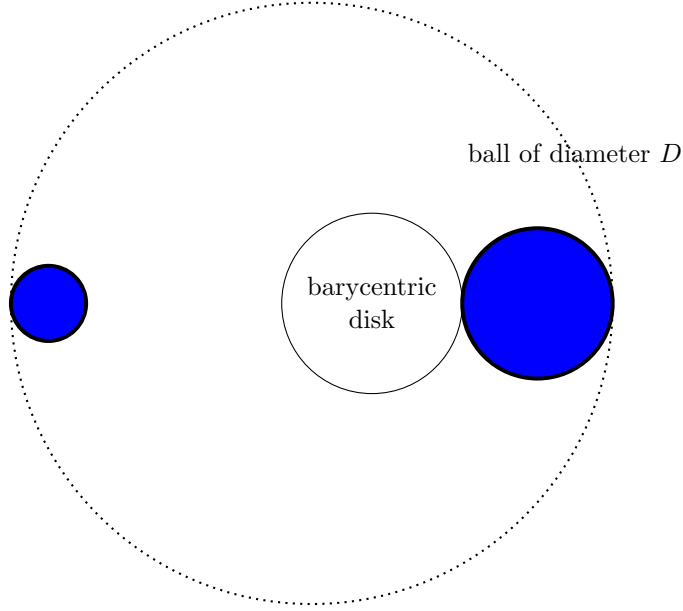
\end{center}
As we will see, this situation occurs for one and only one value of  $R_1$, depending on the diameter $D$. We will denote this value by $R_1^*(D)$
or $R_1^*$.  By the choice we did, we have $J(K)=(R_1^*+\sqrt{1-{R_1^*}^2}-1)/4$.
Since the function $R_1 \mapsto (R_1+\sqrt{1-R^2_1} -1)/4$ is strictly decreasing, we can find the unique value of $R_1$ such that 
$J(K)=\tau^*$: this is $\widehat{R_1} \simeq 0.881075$. In other words, as soon as we have 
$$R_1^*(D) \geq 0.8811 \ \Rightarrow J(K) < \tau^*$$
then we can claim that a minimizing sequence cannot converge to the ball.

Let us now give an estimate of $R_1^*$. The barycenter of $K$ is on the horizontal line and its abscissa is
$$x_G=R_1^2\left(\frac{D}{2} -R_1\right)+R_2^2\left(-\frac{D}{2} +R_2\right)=-R_1^3+DR_1^2-\frac{D}{2}+(1-R_1^2)^{3/2}.$$
The barycentric ball is tangent to the ball of radius $R_1$ if $x_G+1=\frac{D}{2}-2R_1$ or 
\begin{equation}\label{equ1}
R_1^3-DR_1^2-2R_1+D-1=(1-R_1^2)^{3/2}.
\end{equation}
Let us introduce $q_D(R_1):= R_1^3-DR_1^2-2R_1+D-1$ the left-hand side of \eqref{equ1}. The number $R_1^*$ must satisfy
$q_D(R^*_1)\geq 0$. Now it is easy to see that $R_1 \mapsto q_D(R_1)$ is decreasing on $[0,1]$, $q_D(0)=D-1>0$, $q_D(1)=-2<0$.
Thus $q_D$ has a unique root $\bar{R}$ in $[0,1]$ meaning that $R_1^*\leq \bar{R}$ . Moreover,
$$q_D\left(1-\frac{1}{D}\right)=\frac{1}{D}\left(1-\frac{1}{D}\right)\left(\frac{1}{D}-2\right) <0$$
while
$$q_D\left(1-\frac{1}{D}-\frac{1}{D^2}\right)=\frac{1}{D^3}\left(4-\frac{3}{D^2}-\frac{1}{D^3}\right) >0.$$
Therefore, 
$$1-\frac{1}{D}-\frac{1}{D^2} < \bar{R} < 1-\frac{1}{D}.$$
Assuming $R_1^*\leq \bar{R}$, both sides of \eqref{equ1} are nonnegative, so we can square this equality arriving at the
polynomial expression
$$2R_1^6-2DR_1^5+(D^2-7)R_1^4+(6D-2)R_1^3+(7+2D-2D^2)R_1^2-4(D-1)R_1+D^2-2D=0$$
that can be factorized as
$$(R_1+1)^2[2R_1^4-2(D+2)R_1^3+(D^2+4D-1)R_1^2+(4-2D^2)R_1+D^2-2D]=0.$$
Let us introduce the polynomial of degree 4:
$$p_D(R):=2R^4-2(D+2)R^3+(D^2+4D-1)R^2+(4-2D^2)R+D^2-2D$$ 
and we know that $R_1^*$ is the root of this polynomial less than $\bar{R}$.
Now, an elementary analysis proves that $R\mapsto p_D(R)$ is first decreasing, then increasing on $[0,1]$ with $p_D(0)=D^2-2D >0$
and $p_D(1)=1>0$. We recall that we have chosen $D\geq 5$ here. We also have
$$p_D\left(1-\frac{1}{D}-\frac{1}{D^2}\right)=\frac{2}{D^8}+\frac{8}{D^7}+\frac{8}{D^6}-\frac{2}{D^5}-\frac{5}{D^4}-\frac{2}{D^3} <0 \ \mbox{for }
D\geq 5$$
while 
$$p_D\left(1-\frac{1}{D}-\frac{2}{D^2}\right)=\frac{32}{D^8}+\frac{64}{D^7}+\frac{16}{D^6}-\frac{16}{D^5}-\frac{2}{D^4}-\frac{4}{D^3}
+\frac{1}{D^2} > 0 \ \mbox{for } D\geq 5$$
and 
$$p_D\left(1-\frac{1}{D}\right)=\frac{2}{D^4}-\frac{4}{D^3}
+\frac{1}{D^2} > 0 \ \mbox{for } D\geq 5.$$
Therefore, we deduce that
\begin{equation}\label{encadreR1*}
1-\frac{1}{D}-\frac{2}{D^2} < R_1^* < 1-\frac{1}{D}-\frac{1}{D^2} < \bar{R}
\end{equation}
the other root of the equation being between $1-\frac{1}{D}-\frac{1}{D^2} $ and $1-\frac{1}{D}$.
From the first inequality in \eqref{encadreR1*}, we deduce immediately that $R_1^*(D)\geq 0.8814$ as soon as $D\geq 10.1$.
Observe that, for $D> 10$
$$
p_D\left(1-\frac{1}{D}-\frac{1.8}{D^2}\right)=\frac{20.9952}{D^8}+\frac{46.656}{D^7}
+\frac{15.552}{D^6}-\frac{12.816}{D^5}-\frac{3.4}{D^4}-\frac{3.28}{D^3}
+\frac{0.64}{D^2} > 0\,.
$$
In conclusion, for $D\geq 10$, a sequence converging to the unit ball cannot be minimizing
since we have found a better competitor.
\end{enumerate}

We have seen in the previous steps that a minimizing sequence $K_n$ cannot converge to a ball (when $D\geq 10$).
Since such a sequence satisfies
\begin{itemize}
\item $K_n \subset B$, where $B$ is the ball centered at the origin of diameter $D$,
\item $|K_n|=\pi$,
\item $P(K_n) \leq  2\pi \sqrt{2}$ (see \eqref{perimeterbounded}),
\end{itemize}
we infer (by classical compactness from the embedding of $BV(B)$ onto $L^1(B)$, see e.g. \cite{HP}) that there exists a set
$K$ such that the characteristic function of $K_n$ converges to the characteristic function of $K$ and $P(K)\leq \liminf P(K_n)$.
This implies, in particular, that $|K|=\pi$ and the barycenter of $K_n$ converges to the barycenter of $K$. Therefore 
$$\lambda_0(K_n) \rightarrow \lambda_0(K)$$
while 
$$\delta(K) \leq \liminf \delta(K_n)$$
proving that the set $K$ realizes the minimum of $J$.
\end{proof}

\begin{rem}\label{remarkexistence}
According to  \eqref{threshold},  we can assume that $J(K^*)<0.0885<0.1$, if $K^*$ is a minimizer.
\end{rem}
\section{Regularity and optimality conditions}\label{section4}
\subsection{Regularity}
The important term in the functional $J(K)$ being $\delta(K)$ that is related to the perimeter, we can expect that 
the regularity of the minimizer $K^*$ is strongly related to the techniques to prove the regularity for the classical isoperimetric problem. However, we cannot expect a strong regularity
when $K^*$ crosses the barycentric ball and we have also to take into account parts of the boundary of $K^*$ that touch the boundary of the big ball
of diameter $D$ (if any). Therefore, we claim:
\begin{theo}\label{theoregul}
Let $K^*$ denotes a minimizer for the functional $J$ defined in \eqref{defJ}.
\begin{itemize}
\item The boundary of the minimizer $\partial K^*$ is globally $C^{1,1}$.
\item The parts of the boundary that do not touch the barycentric ball are analytic.
\end{itemize}
\end{theo}
\begin{proof}
According to Remark \ref{remarkexistence}, we can assume that $J(K^*)<0.1$.
Since $\lambda_0(K^*) \leq 2$, this implies that
$$
P(K^*) \leq (1+4\cdot 0.1)\cdot 2\pi \leq 9\,.
$$
To prove the regularity of the statement, we follow the same strategy as M. Cicalese and G.P. Leonardi in \cite{CiLeexistence} using the notion of 
{\it strong $\Lambda$-minimizer of the perimeter} (see Section \ref{section_preliminaries}). We note $E=K^*$.

Here we choose  $R=\sqrt{\lambda_0(E)/12\pi} < 1/4$ and $0<r\leq R$.
This implies
\begin{equation}\label{measure_ball}
|B(x,r)|=\pi r^2\leq \pi R^2\leq \frac{\lambda_0(E)}{12} \leq \frac{1}{6}.
\end{equation}
We fix $x\in E=K^*$. We assume for the moment that  $x \notin \partial D$.
So, let us consider a set $F$ of finite perimeter that coincides with $E$ outside the ball $B_r=B(x,r)$. We will denote $P(F,B)$ and $P(E,B)$
the perimeter of $F$ and $E$ respectively inside the ball $B_r$ and we can assume that $P(F,B)\leq P(E,B)$ otherwise there is nothing to prove.  We also
denote by $B^c$ the complement of $B_r$ and we will use the fact that $F$ and $E$ coincide on $B^c$.
\begin{enumerate}
\item
The starting point is the minimality of $E$ that is,
\begin{equation}\label{minimaE}
\frac{\frac{P(E)}{2\pi} -1}{\lambda_0(E)^2} \leq \frac{\frac{P(F)}{2\sqrt{\pi|F|}} -1}{\lambda_0(F)^2}
\end{equation}
(we use here the fact that the ball of same volume as $F$ has a radius equal to $\sqrt{|F|/\pi}$).
The inequality \eqref{minimaE} is equivalent to
\begin{equation}\label{miniE2}
P(E,B) + P(E\cap B^c) \leq 2\pi \frac{\lambda_0(E)^2}{\lambda_0(F)^2}\left\lbrack  \frac{P(F,B) + P(E\cap B^c)}{2\sqrt{\pi |F|}} \,-1\right\rbrack +2\pi.
\end{equation}
First of all, recalling that $|F\Delta E|=|F\setminus E| + |E\setminus F|$ and
$|F|-|E|= |F\setminus E| - |E\setminus F|$, we can write
\begin{equation}\label{areaF}
\pi + |E\Delta F| \geq |F| \geq \pi - |E\Delta F|.
\end{equation}
Using the inequality $(1-u)^{-1/2} \leq 1+u$ when $0\leq u\leq 1/2$, we obtain from \eqref{areaF} with $u=|E\Delta F|/\pi$
\begin{equation}\label{majF}
\frac{1}{2\sqrt{\pi |F|}} \leq \frac{1}{2\pi} \left(1-\frac{|E\Delta F|}{\pi}\right)^{-1/2} \leq  \frac{1}{2\pi} \left(1+\frac{|E\Delta F|}{\pi}\right).
\end{equation}
\item
Now we want to compare $\lambda_0(E)$ and $\lambda_0(F)$.
Let $G_F$ and $G_E$
be the barycenters of $F$ and $E$ respectively and $B_F,B_E$ their respective barycentric balls. 
Without loss of generality, we can assume that $G_E=O$ the origin.
Let 
$\chi_K$ denote the characteristic function of the set $K$. Then
\begin{equation}\label{regularitytriangular}
\pi \lambda_0(E)=\|\chi_E-\chi_{B_E}\|\leq \|\chi_E -\chi_F\| + \|\chi_F - \chi_{B_F}\| + \|\chi_{B_F} - \chi_{B_E}\|.
\end{equation}
Now, $\|\chi_E -\chi_F\|=|E\Delta F|$, $\|\chi_F - \chi_{B_F}\|=\pi\lambda_0(F)$ and it remains to estimate $\|\chi_{B_F} - \chi_{B_E}\|$
by using Lemma \ref{lemmapalle}. For that purpose, we need an estimate of the distance between $G_F$ and $G_E$. We have, since $E=F$ in $B_r^c$
$$
G_F=\frac{1}{|F|} \int_F X dX = \frac{1}{|F|} \,  \left(\int_{E\cap B_r^c} XdX + \int_{F\cap B_r} XdX\right)
$$
and similarly for $G_E$:
$$
O=G_E=\frac{1}{|E|} \int_E X dX = \frac{1}{|E|} \,\left(\int_{E\cap B_r^c} XdX + \int_{E\cap B_r} XdX\right)
$$
which implies that
$$
\int_{E\cap B_r^c} XdX =- \int_{E\cap B_r} XdX
$$
Therefore
$$
G_F=G_F-O=\frac{1}{|F|} \,\left( \int_{F\cap B_r} XdX - \int_{E\cap B_r} XdX \right)
$$
It follows that, denoting $E_1=E\cap B_r$, $F_1=F\cap B_r$:
$$|G_F|=|G_F-G_E| \leq \frac{1}{|F|} 2r \left( |E_1\setminus F_1| + |F_1\setminus E_1| \right) 
= \frac{2r}{|F|} \,|E\Delta F|
\leq 
\frac{4r}{\pi} |E\Delta F|
$$
by using \eqref{areaF}. 

Recall that the balls $B_E$ and $B_F$ have not the same volume. So, let $B'_F$ the ball centered at $G_F$
with area $\pi$. We have
\begin{equation}\label{balls}
|B_E\Delta B_F| \leq |B_E\Delta B'_F|+|B'_F\Delta B_F|=|B_E\Delta B'_F| +||F|-\pi|\leq |B_E\Delta B'_F| +|E\Delta F|
\end{equation}
where we have used \eqref{areaF} for the last inequality.
Now it remains to use Lemma \ref{lemmapalle} to estimate $|B_E\Delta B'_F|$. 
We introduce the function $f: a \mapsto 4a\arcsin(\frac{a}{2})+2a\sqrt{1-\frac{a^2}{4}}$ and we note that $f$
is increasing and $f(a)\leq 4a$ when $a\in [0,1]$. Since $\frac{4r}{\pi} |E\Delta F| \leq \frac{|E\Delta F|}{\pi}<1$ (here we have used $r\leq 1/4$), 
we have
$$
|B_E\Delta B'_F| \leq f\left(\frac{4r}{\pi}|E\Delta F|\right)
\leq 
4\cdot \frac{4r}{\pi}|E\Delta F|
\leq \frac{4}{\pi} |E\Delta F|
$$ 
Therefore, by (\ref{balls}),
 $$
 |B_E\Delta B_F| \leq \left(1+\frac{4}{\pi}\right) |E\Delta F|\leq 3 |E\Delta F|.
 $$
By (\ref{regularitytriangular}),
this implies the estimate
\begin{equation}\label{dernierelambda_0}
\pi\lambda_0(E) \leq  \pi\lambda_0(F) + 4 |E\Delta F|\leq
\pi\lambda_0(F) + 4 \pi r^2
\end{equation}
since $E\Delta F\subset \subset B(x,r)$.
\item
It follows from (\ref{dernierelambda_0}) that 
$$
\frac{\lambda_0(E)}{\lambda_0(F)}\leq 1+\frac{4|E\Delta F|}{\pi\lambda_0(F)} \leq 1+\frac{4|E\Delta F|}{\pi(\lambda_0(E)-4\pi r^2)}\leq 
1+\frac{6 |E\Delta F|}{\pi\lambda_0(E)}
$$
since $\pi r^2\leq \lambda_0(E)/12$ by (\ref{measure_ball}).
Squaring the previous inequality, and using the fact that $(1+u)^2\leq 1+3u$ as soon as $0\leq u\leq 1$, we infer
\begin{equation}\label{ratiolambda0}
\left(\frac{\lambda_0(E)}{\lambda_0(F)}\right)^2 \leq 1+\frac{18 |E\Delta F|}{\pi\lambda_0(E)} .
\end{equation}
Finally, replacing in \eqref{miniE2}, the different estimates obtained in  \eqref{majF} and \eqref{ratiolambda0} 
and using also $P(F) \leq P(E), |E\Delta F|\leq \pi$, we obtain
$$
P(E,B)\leq
2\pi+
2\pi\left(1+\frac{18|E\Delta F|}{\pi \lambda_0(E)}\right)\left[
-1+\frac{P(F,B)}{2\pi}\left(1+\frac{|E\Delta F|}{\pi}\right)
\right]
$$
and therefore
$$
P(E,B) \leq P(F,B) + \Lambda |E\Delta F|
$$
with
$$\Lambda=\frac{36}{\pi\lambda_0(E)} \left(P(E) -\pi\right) + \frac{P(E)}{\pi}.$$
This proves the $C^{1,1}$ regularity of  the minimizer $K^*$ outside $\partial D$, as explained in Section \ref{section_preliminaries}.
\end{enumerate}
Even if the boundary of $K^*$ meets the boundary of the ball of diameter $D$,
the global $C^{1,1}$ regularity follows by Theorem \ref{thm_stred-zi}.
The analyticity of the parts of the boundary that do not cross the boundary of the barycentric ball still follows from Theorem \ref{thm_stred-zi}.
\end{proof}
\begin{rem}
The analyticity also follows from the next section where we will be able to prove that the curvature $\mathcal{C}$
(that is defined almost everywhere since $\partial K^*$ is $C^{1,1}$) satisfies the optimality condition:
$$\mathcal{C}=C_0 +\mu_1 x +\mu_2 y$$
where $C_0$ is a constant (different inside and outside the barycentric ball) and $\mu_1,\mu_2$ two kind of Lagrange multipliers associated to
constraints on the barycenter. From this relation we see that $\mathcal{C}$ is not only $C^\infty$ but also analytic. The analyticity of the
boundary follows.
\end{rem}

\subsection{Optimality conditions}
We are going to write the optimality conditions on the boundary of the optimal set $K^*$. For that purpose, we use the
framework of shape derivative as explained for example in \cite[chapter 5]{HP}. We essentially reproduce here the computations done in our 
previous paper \cite{BCH_Annali}.
We go on assuming that all the considered sets have area equal to $\pi$. 

\begin{theo}\label{ocPrinceton}
Let $K$ be an optimal set minimizing the functional $J$ and assume its barycenter is at the origin. Let $B$ be the unit ball centered at the origin. Let $\partial K^{IN}=\partial K \cap B$, $\partial K^{OUT}=\partial K \cap B^c$, $\partial B^{IN}=\partial B \cap K$, 
$\partial B^{OUT}=\partial B \cap K^c$. \\
Then at every boundary point $(x,y)$ of $\partial K$that are not on the boundary of $B$ and not on the boundary of the ball 
of diameter $D$ containing $K$, the curvature $C(x,y)$ satisfies:
\begin{equation}\label{opcond}
C(x,y)=1-3\delta(K)+\frac{4\delta(K)}{2\pi \lambda_0(K)} \left(|\partial B^{OUT}|-|\partial B^{IN}|\right)
\pm \frac{4\delta(K)}{\lambda_0(K)} +\mu_1 x + \mu_2 y\,, 
\end{equation} 
($+$ at the exterior of $B$ and $-$ in the interior of $B$)
where
\begin{equation}\label{defmu1}
\mu_1= \frac{4\delta(K)}{\pi\lambda_0(K)} \left[\int_{\partial B^{OUT}} \cos t dt -  \int_{\partial B^{IN}} \cos t dt \right]\,,
\end{equation}
\begin{equation}\label{defmu2}
\mu_2=\frac{4\delta(K)}{\pi\lambda_0(K)} \left[\int_{\partial B^{OUT}} \sin t dt -  \int_{\partial B^{IN}} \sin t dt \right]\,.
\end{equation}
\end{theo}

\begin{proof}
Let us remark that, according to the regularity results stated in Theorem \ref{theoregul}, the curvature is defined almost everywhere (and bounded)
and the following computations are well justified.
For simplicity we denote now $\delta$ and $\lambda_0$ for $\delta(K)$ and $\lambda_0(K)$.

 Let $V$ be a perturbation, that is, a smooth map $V:\mathbb{R}^2\to\mathbb{R}^2$. We denote by $V\cdot n$ the scalar product of $V$ with the outer unit normal vector $n$ to $\partial K$.
Let $K_t=(I+t V)(K)$. Then the area of $K_t$ satisfies
$$
|K_t|=\pi +t \int_{\partial K}  V\cdot n +o(t)\,.
$$
Now we look at the barycenter $G_t=(x_t,y_t)$ of $K_t$.
We have:
 $$\displaystyle \int_{K_t} x\;dxdy= 0 +t\int_{\partial K} x V\cdot n +o(t).$$
Since by definition
$\displaystyle 
x_t=\frac{1}{|K_t|}\int_{K_t} x\; dxdy
$,
by the above formulas one has
$$
x_t=\frac{t}{\pi}\int_{\partial K} x\;V\cdot n  +o(t)\,.
$$
A similar formula holds for $y_t$:
$$
y_t=\frac{t}{\pi}\int_{\partial K} y\;V\cdot n +o(t)\,.
$$
Now, let  $B_t=(I+tW)(B)$, where
$$
\displaystyle W(x,y)=(a,b)+\alpha (x,y)\,,
$$
with
$$
(a,b)=\frac{1}{\pi}\left(\int_{\partial K} x V\cdot n, \int_{\partial K} y V\cdot n \right)
\,, \qquad
\alpha = \frac{1}{2\pi} \int_{\partial K} V\cdot n\,.
$$
The ball $B_t$ is, at the first order, the barycentric ball of $K_t$.
The difference between  $|K_t \Delta B_t|$ and $|K \Delta B|$ is given by two terms:
$$\displaystyle 
|K_t \Delta B_t|-|K \Delta B|= \pm \; t \int_{\partial B} W\cdot n \pm t \int_{\partial K} V\cdot n\,:
$$
for the first term of the right hand side  $+$ is on  $\partial B^{OUT}$ and 
$-$ is on $\partial B^{IN}$;
for the last term of the right hand side,  $+$ is $\partial K^{OUT}$ and 
$-$ on $\partial K^{IN}$.

In the next part of the proof we will write
$\displaystyle 
|K_t \Delta B_t|=|K \Delta B|+  t R
$.
\item
We have
$$\lambda_0(K_t)=\frac{|K_t \Delta B_t|}{|K_t|}=\frac{|K \Delta B| + tR +o(t)}{|K| + t\int_{\partial K} V\cdot n +o(t)}
$$
$$
= \lambda_0(K)\cdot \frac{1+t\frac{R}{|K \Delta B|}}{1+ \frac{t}{\pi}\int_{\partial K} V\cdot n} +o(t)= \lambda_0(K) + t\left[\frac{R}{\pi}-\frac{\lambda_0}{\pi} \int_{\partial K} V\cdot n\right]\,+o(t).
$$
Therefore, the first derivative of $t\mapsto \lambda_0(K_t)$ is :
$$
d\lambda_0(K,V)=\frac{1}{\pi} \left[\pm \int_{\partial B} W\cdot n \pm \int_{\partial K} V\cdot n - \lambda_0 \int_{\partial K} V\cdot n   \right]\,,
$$
that is,
\begin{eqnarray*}
d\lambda_0(K,V)=\frac{1}{\pi} \left[ \int_{\partial B^{OUT}} W\cdot n - \int_{\partial B^{IN}} W\cdot n\right]+ \\
 +\frac{1}{\pi}  \left[\int_{\partial K^{OUT}} V\cdot n - \int_{\partial K^{IN}} V\cdot n - \lambda_0 \int_{\partial K} V\cdot n   \right]\,.
\end{eqnarray*}
If $r_t$ is the radius of the ball having the same area as  $K_t$, then 
$$
r_t=\sqrt{\frac{\pi + t\int_{\partial K} V\cdot n}{\pi}} +o(t)=1+t\; \frac{\int_{\partial K} V\cdot n}{2\pi} +o(t).
$$ 
This gives
$$
\delta(K_t)=\frac{P(K_t)}{2\pi r_t} -1=\frac{P(K_t)}{2\pi + t \int_{\partial K} V\cdot n+o(t)} -1=\frac{P(K) + t\int_{\partial K} \mathcal{C}\ V\cdot n}{2\pi + t \int_{\partial K} V\cdot n} -1 +o(t)
$$
where $\mathcal{C}$ denotes the curvature of the boundary.
With the same computations as for $\lambda_0$
$$
\delta(K_t)=\delta(K) + t\frac{\int_{\partial K} \mathcal{C}\ V\cdot n}{2\pi} -t\frac{P(K) \int_{\partial K}V\cdot n}{4\pi^2}  +o(t)
$$
and so, the first derivative of $t\mapsto \delta(K_t)$ is given by
$$
d\delta(K,V)=\int_{\partial K} \left[\frac{\mathcal{C}}{2\pi}-\frac{(\delta +1)2\pi}{4\pi^2} \right]V\cdot n = 
\int_{\partial K} \frac{\mathcal{C}- \delta -1}{2\pi} V\cdot n\,.
$$
The optimality condition for the functional $J$ 
\begin{equation}\label{opcond2}
\displaystyle
\frac{d\delta}{\lambda_0^2}-\frac{2\delta}{\lambda_0^3} d\lambda_0=0,
\end{equation}
can be written as
$$\int_{\partial K} (\mathcal{C}-\delta - 1)V\cdot n = \frac{4\delta}{\lambda_0} \left [\int_{\partial B^{OUT}} W\cdot n - \int_{\partial B^{IN}} W\cdot n\right] 
$$
$$+
 \frac{4\delta}{\lambda_0} \left[\int_{\partial K^{OUT}} V\cdot n - \int_{\partial K^{IN}} V\cdot n 
- \lambda_0 \int_{\partial K} V\cdot n \right] \,.
$$
Now,
$W\cdot n=a\cos\theta + b\sin\theta +\alpha$ (since $(x,y)\cdot n=1$ on $\partial B$), so
$$
W\cdot n =\cos\theta \int_{\partial K} x V\cdot n +\sin\theta \int_{\partial K} y V\cdot n
+\frac{1}{2\pi} \int_{\partial K} V\cdot n \,.
$$
Therefore, since the identity \eqref{opcond2} must hold for any deformation field $V$, we finally get
\begin{equation}\label{optimcond}
\mathcal{C}=\delta+1-4\delta+\frac{4\delta}{2\pi \lambda_0} \left(|\partial B^{OUT}|-|\partial B^{IN}|\right)
\pm \frac{4\delta}{\lambda_0} +\mu_1 x + \mu_2 y
\end{equation}
which is the expected result.
\end{proof}
\begin{rem}\label{pendulum}
By a rotation, we can assume that the Lagrange multiplier $\mu_2=0$.
Denoting by $\theta$ the angle of the tangent of the boundary of the optimal set
with the horizontal
axis and by $s$ the curvilinear abscissa, the optimality condition \eqref{opcond} can be written
as
$$
\frac{d \theta}{ds} = a +\mu_1 x(s)
$$
with a constant $a$ different inside or outside the barycentric disk.  Differentiating this relation with respect to  $s$ yields the pendulum equation
\begin{equation}\label{eq_pendule}
\frac{d^2 \theta}{ds^2} = \mu_1 \cos\theta(s).
\end{equation}
\end{rem}
\section{Qualitative properties}\label{section5}
In this section, we give some properties of the minimizer $K^*$. Let us remark first that up to a rotation, we can always assume that the Lagrange multiplier $\mu_2$ defined in \eqref{defmu2} is equal to zero.
 We denote by $B_D$ the disk of diameter $D$ that contains our optimal domain $K^*$
and by $B_G$ the barycentric disk of $K^*$.
\subsection{Connected components}
 We will say that two (or more) connected components touch the boundary of $B_D$
in {\it opposite} points when it is not possible to move the disk $B_D$ without making at least one component going outside $B_D$.
Let us start with this preliminary remark:
\begin{rem}\label{diammax}
Each connected component $K_i$ of $K^*$ has a diameter bounded:
\begin{equation}\label{bornediam}
diam(K_i) \leq 4.26\,.
\end{equation}
Indeed, recalling that by (\ref{threshold}) $\delta(K^*)/\lambda_0^2(K^*)\leq 0.089$, and $\lambda_0(K^*) \leq 2$, we have
$\delta(K^*)\leq 4\cdot 0.089=0.356$ and therefore $P(K_i)\leq P(K^*)\leq 0.356\cdot 2\pi+2\pi \leq 8.52$. Now, $K_i$ being connected, denoting by $conv(K_i)$
its convex hull, we have
$$diam(K_i)=diam(conv(K_i))\leq \frac{P(conv(K_i))}{2} \leq \frac{P(K_i)}{2} \leq 4.26.$$
\end{rem}
\begin{theo}\label{theotwocomp}
There exists a number $D^{**}$ such that for $D\geq D^{**}$, the minimizer $K^*$ has exactly two connected components.
Moreover these two connected components touch the boundary of $B_D$
in {\it opposite} points.
\end{theo}
\begin{proof}
{\bf (1)}
Let us remark that the numerical sequence
$$
m_D:=\min\{J(K), diam(K)\leq D\}
$$
is non increasing with respect to $D$ (since the class of competitors enlarges, as $D$ increases).
Assume there exists a sequence $D_n$ going to $+\infty$
such that a minimizer $K_n^*$, with the diameter $D_n$ has only one connected component.
Now, the bound \eqref{bornediam}, with the fact that $D_n\geq 10$, shows that $K_n^*$ cannot touch the disk $B_{D_n}$
 in points located in two different half-planes
(passing through the center of $B_{D_n}$). In other words, it is possible to move inside $B_{D_n}$ the set $K_n^*$ in such a way that it does not touch 
$B_{D_n}$. But that would mean that $K_n^*$ is necessarily fixed and in particular, the sequence $m_{D_n}$ would be constant.
Since this sequence is non increasing, this would imply that $m_D$ itself is constant.
This is in contradiction with  Remark \ref{limite_Fuglede}.

\medskip\noindent
{\bf (2)}
In the sequel, we now assume that $D\geq D^{**}$.
Now, let us assume that the optimal set $K^*$ has three connected components $K_0,K_1,K_2$
(we will treat the case of more than three connected components at the
end of this proof).
Each connected component can be in three different situations: outside the barycentric disk $B_G$, crossing the  barycentric disk or inside the
 barycentric disk.  By "crossing the barycentric disk", we mean that both intersections of the component with $B_G$ and with its complement 
 have positive measure. We have also to take into account the cases where one (or more) component touches the boundary of the big disk $B_D$.
The strategy of the proof (by contradiction) is to be able to be reduced to one of the following situations:
\begin{description}
\item[A] two components that are "free" outside $B_G$ (by free we mean that they do not touch $\partial B_G$ and $\partial B_D$,
\item[B] two components that are "free" inside $B_G$ (that is, they do not touch $\partial B_G$)
\item[C] one component that crosses $B_G$ that we can translate outward $B_G$ without changing the barycenter.
\end{description}
Indeed, in the case A, if $K_1,K_2$ are these two components, without loss of generality, we can assume that
$$\frac{|K_1|}{P(K_1)} \frac{P(K_2)}{|K_2|} \leq 1.$$
Then, we can expand $K_2$ with a factor $t_2>1$ and shrink $K_1$ with a factor $t_1<1$ such that the total volume is preserved, i.e.
$$(t_1^2-1)|K_1|+(t_2^2-1)|K_2|=0.$$
Up to some translation, we can assume that the barycenter of $t_1K_1 \cup t_2 K_2$ is the same as $K_1 \cup  K_2$, therefore $\lambda_0$ is unchanged. Finally, the
new perimeter being $t_2P(K_2)+t_1P(K_1)$ is strictly smaller than $P(K_2)+P(K_1)$ because
$$(t_2-1)P(K_2)<(1-t_1)P(K_1) \ \Leftrightarrow \ \frac{t_1+1}{t_2+1 } \frac{|K_1|}{P(K_1)} \frac{P(K_2)}{|K_2|} < 1$$
giving the desired contradiction.

In the case B, we can simply replace the two components by a disk of same volume (and same barycenter), keeping $\lambda_0$ unchanged and
decreasing the perimeter. Even if this new disk is not completely contained in $B_G$ this is favourable since, in that case, $\lambda_0$ would
increase.

Finally, in the case C, translating the crossing component without changing the barycenter increases $\lambda_0$, giving the contradiction.

Considering the three components, here are all the ten possibilities:
 \begin{enumerate}
 \item $K_0$ crosses $B_G$, $K_1$ is outside $B_G$, $K_2$ is inside $B_G$,
 \item $K_0$ crosses $B_G$,  $K_1$ and $K_2$ are outside $B_G$,
 \item $K_0$ crosses $B_G$,  $K_1$ and $K_2$ are inside $B_G$,
 \item $K_0$ and $K_1$ cross $B_G$, $K_2$ is outside $B_G$,
 \item $K_0$ and $K_1$ cross $B_G$, $K_2$ is inside $B_G$,
 \item $K_0$ and $K_1$ are inside $B_G$, $K_2$ is outside $B_G$, 
 \item $K_0$ and $K_1$ are outside $B_G$, $K_2$ is inside $B_G$,
\item $K_0,K_1$ and $K_2$ are outside $B_G$,
\item $K_0,K_1$ and $K_2$ are inside $B_G$,
 \item $K_0,K_1$ and $K_2$ cross $B_G$.
\end{enumerate} 
The case (i) is impossible since the total area of $K^*$ is $\pi$.\\
The cases (c) and (f) can be treated with case [B].\\
In the cases (e) and (j), by the diameter bound \eqref{bornediam} it is impossible that two components touch $B_D$ in points located in two different half-planes
(passing through the center of $B_D$). Therefore, we can move $K^*$
inside $B_D$ and then we are led to the favourable situation [C] by choosing two crossing components and moving them outward.\\
In the case (a), the two components $K_0$ and $K_1$ must touch the boun\-da\-ry of the disk $B_D$ in two "opposite" points
otherwise they can be moved outward. Now we can use the third component $K_2$ 
(possibly after a rotation of $K_0,K_2$ if needed)
 and we move it together with $K_1$ in such a way that they are closer to $K_0$, infinitesimaly. In that process, either $K_2$ goes out of $B_G$ and we have increased
 $\lambda_0$ or $K_2$ stays inside $B_G$, but now we have room to move $K_0$ outside $B_G$ while $K_1$ moves in the other way, increasing
 one more time $\lambda_0$.\\
 In the cases (b) and (d), we can assume that  $K_1$ and $K_2$ touch the boundary of $B_D$ (and possibly also $K_0$) in opposite points. 
 In that case, assuming that it is $K_1$ that is close to $K_0$ we start to rotate $K_0$ and $K_1$ by a rotation of center $O$
 (the center of $B_D$,)  in two opposite directions. In that case their barycenter moves in the direction of $K_2$ and since $K_0$ and $K_1$
 remain stuck to the boundary of $B_D$ the quantity $\lambda_0$ increases.  \\
 Finally, in the cases (g) and (h), we use a similar idea: first a rotation along the boundary of $B_D$ of two components $K_1,K_2$. Then, since the 
 barycentric disk $B_G$
has moved letting free $K_1$ we move $K_0$ and $K_1$ closer to each other, making the strategy [A] possible.

\medskip\noindent
{\bf (3)}
Now, let us assume that we have more than three connected components $K_0,K_1,K_2,K_3,\ldots$. As soon as we can use one of the
strategy described before with three components, we are done. 
The only action that would not be possible with more components would be to move $B_D$ because a supplementary component may prevent to do this.
But examining the different strategies presented in the case of three components, we can see that we never use such a trick.
Therefore, we are done also in that case.
 
Finally the claim about the fact that the two connected components must touch the boundary of the ball $B_D$ in two opposite points
comes from the following fact: if it was not the case, by moving $B_D$ we could  make $K_0$ and $K_1$ enter into the interior of $B_D$.
Then we have room to translate both components and apply one of the strategy [A] or [C] to reach a contradiction.

\end{proof}

\subsection{Nature of the boundary}
It is generally difficult to identify the optimal domain for quantitative inequalities. The only example we have in mind
is the "stadium" found by S. Campi in \cite{Ca} and Alvino, Ferone, Nitsch in \cite{AFN} for the case of convex domains
and the Fraenkel asymmetry in the plane.  It turns out that, in most cases (proved or conjectured), 
the optimal domain (in the plane) has a boundary composed of different arcs of circle (including possibly segments).
Therefore, the following theorem is quite surprising.
\begin{theo}
The boundary of the minimizer does not contain any arc of circle. In other terms, the Lagrange multiplier $\mu_1$ defined in \eqref{defmu1} is not zero.
\end{theo}
\begin{proof}
In view of the optimality conditions \eqref{optimcond} expressing the curvature, we see that the presence of an arc of circle (either inside
or outside the barycentric ball) is equivalent to $\mu_1=\mu_2=0$. As soon as there is an arc of circle somewhere, all parts of
the boundary of $K^*$ are composed of arcs of circle, with the same radius, say $R_0$ outside and the same radius $R_1$ inside.
In the proof, we argue by contradiction.

\medskip\noindent
According to Theorem \ref{theotwocomp}, if we assume $D\geq D^{**}$, the only case we have to consider is the case of two components $K_0,K_1$.

Each connected component can be in three different situations: outside the barycentric disk $B_G$, crossing the  barycentric disk or inside the
 barycentric disk.  Moreover, we know that both $K_0$ and $K_1$ must touch the boundary of the disk $B_D$ in opposite points, by Theorem \ref{theotwocomp}.
 Therefore, the only situations we have to consider are
\begin{enumerate}
\item $K_0$ is outside $B_G$ and $K_1$ is inside $B_G$,
\item $K_0$ is outside $B_G$ and $K_1$ crosses $B_G$,
\item $K_0$ and $K_1$ are outside $B_G$.
\end{enumerate}
Let us first consider the component $K_0 \subset B_D\setminus B_G$. If its contact with $\partial B_D$ has only one point then $K_0$
must be a disk. If not, this contact is an arc of $\partial B_D$ and the remaining boundary of $K_0$ is another arc of circle $\gamma_0$.  But since
$K_0$ has global regularity $C^1$, the only possibility would be that $\gamma_0$ has the same curvature than $\partial B_D$: this is impossible.
Therefore we can conclude that $K_0$ is a disk touching $\partial B_D$ in only one point.

In case 3, we can argue in the same way for $K_1$ and therefore $K_1$ is also a disk with the same radius $1/\sqrt{2}$ according to the optimality condition
and the area constraint.
In that case, $\delta(K^*)/\lambda_0^2(K^*)=(\sqrt{2}-1)/4$ being a non competitive value, by (\ref{threshold}).

In the case 1,  both domains are disks of respective radii $R_0$ and $R_1$. Since $K_1$ is inside the barycentric disk, necessarily $R_1\geq 1/\sqrt{2}\geq R_0$.
Now an explicit computation gives
$$\delta(K^*)=R_0+R_1-1=R_1+\sqrt{1-R_1^2}-1 , \quad \lambda_0(K^*)=2(1-R_1^2)$$
and a simple study of the function $R_1\mapsto \delta/\lambda_0^2$ shows that the minimum is achieved by $R_1=1/\sqrt{2}$ and is equal to
$\sqrt{2}-1$ that is not competitive.

The case 2 is more complicated. Assume first that $K_1$ touches the boundary of $B_D$ in only one point. As already observed  in the proof of Theorem 4.4 of 
\cite{BCH_COCV},  the succession of arcs of circle outside the barycentric disk with curvature $\mathcal{C}^{OUT}$ and arcs of circle
inside the barycentric disk, with curvature $\mathcal{C}^{IN}$, the junction of these arcs being $C^1$ is only possible with a periodicity i.e. a repetition
of $N$ succession of such a couple of arcs.  In that case, this component $K_1$ has a rotational symmetry with respect to the center $G$ of $B_G$
meaning that $G$ is the barycenter of $K_1$: thus it cannot be the barycenter of $K^*=K_0\cup K_1$: a contradiction.
Now assume that $K_1$ touches the boundary of $B_D$ along an arc of circle.  Without loss of generality, assume that this arc of circle is centered at the
point $(D/2,0)$.
The remaining part of $\partial K_1$ still satisfies the geometric property
just described. In other words, if $\widehat{K_1}$ denotes the same rotationally symmetric set as in the previous case, $K_1$ is simply the intersection
of $\widehat{K_1}$ with $B_D$. Consequently the barycenter of $K_1$ is on the left of $G$ (it has a lower abscissa). But
since the other component $K_0$ touches $\partial B_D$ in an opposite point, the barycenter of $K^*=K_0\cup K_1$ cannot be equal to $G$: a contradiction.
 
\end{proof}

\begin{rem}
Assuming, for a contradiction, that the boundary of $K^*$ is composed of arcs of circle,
we can prove in a different way that the minimizer cannot have only one connected component for any value of
the diameter $D$ (not only for $D$ large enough), with a similar proof to that one used in \cite{BCH_COCV}.
\end{rem}

\section{Perspectives}
It would be interesting to have a more precise description of the minimizer. For example, the following questions may be considered:
\begin{description}
\item[Symmetry] as observed in Remark \ref{pendulum},
since  the Lagrange multiplier $\mu_2$ can be assumed to be 0, it is reasonable to think
that $K^*$ is symmetric with respect to the horizontal axis. Here are some remarks about a proof of the symmetry, based on the  pendulum equation (\ref{eq_pendule}), that is, 
$$\frac{d^2 \theta}{ds^2} = \mu_1 \cos\theta(s).$$
Starting from a point $A_0=(x_0,y_0)$ where the tangent is vertical ($\theta(0)=\pi/2$),  we can prove, using uniqueness
of the solution, that $\theta(-h)=\pi - \theta(h)$. This implies that the curve is symmetric with respect to the line $y=y_0$
in a neighborhood of $A_0$. To be able to conclude to the full symmetry, it remains to prove that $A_0$ is itself on the
horizontal axis.
\item[Contact with the ball $B_D$] We think that the optimal set $K^*$ touches the ball $B_D$ not only on one point
but along some arc on each side. 
\item[Position of the barycentric disk] The barycentric disk has necessarily a not empty intersection with one component.
Indeed, if it was not the case, we can always replace the two components with two balls with the same volume
keeping possibly after some translation the barycentric disk outside the two components. By this transformation
we keep $\lambda_0$ unchanged, while strictly decreasing $\delta$. Moreover, from the bound on the diameter \eqref{bornediam} 
and the fact that each connected component touches the boundary of $B_D$,
we see that only one component can touch the barycentric disk (at least for $D\geq 4.26+2+4.26=10.52$).
Now, it remains to know how this component, say $K_1$, crosses the barycentric disk $B_G$. We expect that the intersection
$K_1 \cap \partial B_G$ contains only two points.
\end{description}

\section*{Acknowledgements}
The research of  G. Croce and A. Henrot has been supported by the two ANR Projects SHAPO and STOIQUES.


\begin{thebibliography}{10}


\bibitem{AFN} {A. Alvino, V. Ferone, C. Nitsch},
A sharp isoperimetric inequality in the plane. J. Eur. Math. Soc. (JEMS)
13 (2011), 185-206.


\bibitem{BCH_COCV}{C. Bianchini, G. Croce, A. Henrot},
On the quantitative isoperimetric inequality in the plane.
ESAIM: COCV 23 (2017), 517-549.


\bibitem{BCH_Annali}{C. Bianchini, G. Croce, A. Henrot},
On the quantitative isoperimetric inequality in the plane with the barycentric asymmetry,
Ann. Sc. Norm. Super. Pisa, Cl. Sci. 24 (2023), 2477-2500.



\bibitem{Ca}{S. Campi},
Isoperimetric deficit and convex plane sets of maximum translative discrepancy.
Geom. Dedicata 43 (1992), 71-81.


\bibitem{CiLe} {M. Cicalese, G. P. Leonardi}, {A selection principle
for the sharp quantitative isoperimetric inequality}.
Arch. Ration. Mech. Anal. 206 (2012),  617-643.

\bibitem{CiLeexistence} M. Cicalese, G. P. Leonardi,
{Best constants for the isoperimetric inequality in quantitative form}.
J. Eur. Math. Soc. 15 (2013),  1101-1129.

\bibitem{DL}
M. Dambrine, J. Lamboley, Stability in shape optimization with second variation,
Journal of Differential Equations 267, 5 (2019), 3009-3045.



\bibitem{Ferriero-Fusco}
A. Ferriero, N. Fusco,
A note on the convex hull of sets of finite perimeter in the plane, Discrete Cont. Dyn. Syst. Ser. B 11 (2009), 103-108.


\bibitem{FiMP}  
A. Figalli, F. Maggi, A. Pratelli, 
{A mass transportation approach to quantitative isoperimetric inequalities}.
Invent. Math. 182 (2010), 167-211.


\bibitem{Fu89Transactions}{B. Fuglede},
Stability in the isoperimetric problem for convex or nearly spherical domains in $\R^n$. Trans.
Amer. Math. Soc. 314 (1989), 619-638.



\bibitem{Fu93Geometriae}{B. Fuglede},
Lower estimate of the isoperimetric deficit of convex domains in $\R^n$ in terms of asymmetry.
Geom. Dedicata 47 (1993), 41-48.



\bibitem{Fuscopreprint}
N. Fusco, The quantitative isoperimetric inequality and related topics. Bull.  Math. Sciences 5 
(2015), 517-607.

\bibitem{FGP}
N. Fusco, M. S. Gelli, G. Pisante,
On a Bonnesen type inequality involving the spherical deviation.
J.  de Math. Pures et Appliqu\'ees
98 (2012), 616-632.



\bibitem{FuMP} 
{N. Fusco, F. Maggi,  A. Pratelli},
{The sharp quantitative isoperimetric inequality}. 
Ann. of Math.  168
(2008),  941-980.

\bibitem{GP24} {C. Gambicchia, A. Pratelli}, {The sharp quantitative barycentric isoperimetric inequality for bounded sets},
Annali Scuola Normale Superiore - Classe de Scienze  doi.org/$10.2422/2036-2145.202407_013$.


\bibitem{HHW}
R.R. Hall, R. R., 
W.K. Hayman, W. K.,  
A.W. Weitsman, On asymmetry and capacity.
J. Anal. Math. 56 (1991), 87-123.

\bibitem{H}
R.R. Hall, A quantitative isoperimetric inequality in n-dimensional space. 
J. Reine Angew. Math. 428
(1992), 161-176.




\bibitem{HP} {A. Henrot, M. Pierre},   Shape variation and optimization. A geometrical analysis.
EMS Tracts in Mathematics 28. Z\"urich: European Mathematical Society.



\bibitem{LiCRAS} 
G. Li,  X. Zhao, Z. Ding, R. Jiang,
An analytic proof of the planar quantitative isoperimetric inequality.
C. R. Math. Acad. Sci. Paris 353 (2015), 589-593.


\bibitem{Maggi} {F. Maggi,}  Sets of Finite Perimeter and Geometric Variational Problems: An Introduction to Geometric Measure Theory, Cambridge studies in Advanced Mathematics, 135, Cambridge (2012).



%
\bibitem{stred-zi} {E. Stredulinsky and W. P. Ziemer}, Area minimizing sets subject to a volume constraint in a convex set,
J. Geom. Anal. 7 (1997), 653-677. 

\bibitem{tamanini}
{I. Tamanini},
Boundaries of Caccioppoli sets with H\"older continuous normal vector. J. Reine Angew Math. {334} (1982), 27-39.



        




























\end{thebibliography}
\end{document}